\documentclass[a4paper]{article}
\usepackage{ifthen}
\usepackage{pbox}
\usepackage{amsmath, amssymb, amsthm}
\usepackage{pxfonts}
\renewcommand{\coloneq}{\coloneqq}

\newcommand{\bigtimes}{\varprod}
\usepackage[matrix, arrow]{xy}
\usepackage{nameref}
\usepackage[bookmarksnumbered,bookmarksopen,pdfusetitle]{hyperref}

\hypersetup{pdfsubject={characteristic classes},
  pdfkeywords={multiple point manifold; characteristic class; generating
    function}}

\numberwithin{equation}{section}

\newtheorem{theorem}{Theorem}

\newtheorem{corollary}[theorem]{Corollary}

\theoremstyle{definition}
\newtheorem{definition}[theorem]{Definition}

\author{G\'abor Braun\\
  Alfr\'ed R\'enyi Institute of Mathematics\\
  Hungarian Academy of Sciences\\
  Budapest\\
  Re\'altanoda u 13--15\\
  1053\\
  Hungary}

\title{The cobordism class of the multiple points of
  immersions\thanks{Partially supported by Hungarian National Research Fund,
    grant No. T 042 769.}}

\date{}

\newenvironment{diagram}{\begin{lrbox}{\currentdiagram}}{\end{lrbox}\begin{center}
    \usebox{\currentdiagram}
  \end{center}}
\newsavebox{\currentdiagram}

\DeclareMathOperator{\Eq}{Eq}
\DeclareMathOperator{\Deltaalpha}{\Delta\sp{\alpha}}% "diagonal" by alpha
\DeclareMathOperator{\DMalpha}{\Gamma\sp{\alpha}}
\newcommand{\diagonal}[1]{\ensuremath{\Delta\sb{#1}}}

\newcommand{\constant}[2][]{\ensuremath{A_{#2}%shorthand for (-1)^{#2-1}(#2 - 1)!
    \ifthenelse{\equal{#1}{}}%
    {}% Nothing when #1 is empty
    {{#1}^{#2-1}}}}

\newcommand*{\MULTIPLE}[2]{\ensuremath{\widetilde{\Delta}_{#1}(#2)}}% tilde version
\newcommand*{\multiple}[2]{\ensuremath{\Delta_{#1}(#2)}}% in the image

\newcommand{\setZ}{\mathbb{Z}}
\newcommand{\setR}{\mathbb{R}}

\newcommand{\size}[1]{\ensuremath{\lvert#1\rvert}}

\providecommand*{\middle}[1]{\,\vrule\,}

\providecommand*{\hyperref}[2][]{#2}
\providecommand{\coloneq}{\mathrel{:=}}

\providecommand{\bigtimes}{\mathop{\times}}

\begin{document}
\maketitle

\begin{abstract}
  Using generating functions, we derive a multiple point formula for every
  generic immersion \(f\colon M^{m} \looparrowright N^{n}\) between even
  dimensional oriented manifolds.  This produces explicit formulas for the
  signature and Pontrjagin numbers of the multiple point manifolds.  The
  formulas take a particular simple form in many special cases, e.g.\ when
  \(f\) is nullhomotopic, we recover Sz\H ucs's formulas in \cite{MR98i:57058}.
  They also include Hirzebruch's virtual signature formula
  \cite[9.3(4')]{MR34:2573}.
\end{abstract}

\section{Introduction}
\label{sec:introduction}

%% Goal: fat k-fold intersection
Let \(f\colon M^m \looparrowright N^n\) be a generic immersion between oriented
compact even dimensional smooth manifolds.

\paragraph{Aims}
\label{sec:aims}

Let \(M^{\times k}\) be the cartesian product of \(k\) copies of \(M\) and \(f^{\times
  k}\colon M^{\times k} \to N^{\times k}\) be the map induced by \(f\) between the products.
Let
\begin{equation}
  \label{eq:63}
  \Delta \coloneq \{ (x, x, \dotsc , x) : x \in N \} \subseteq N^{\times k}
\end{equation}
be the (narrow) diagonal in \(N^{\times k}\).
Then the \(k\)-tuples of the multiple points of \(f\) form the manifold
\begin{equation}
  \label{eq:1}
  \MULTIPLE{k}{f} \coloneq \{ (x_1, \dotsc , x_k) \in (f^{\times k})^{-1}(\Delta) : x_i \neq
  x_j \text{ if } i \neq j \}.
\end{equation}

Note that by permuting the coordinates, the symmetric group \(S_k\) acts on
\(\MULTIPLE{k}{f}\).  Factoring out with the action leads to the
\emph{manifold of \(k\)-tuple points} of \(f\):
\begin{equation}
  \label{eq:2}
  \multiple{k}{f} \coloneq \MULTIPLE{k}{f} / S_k.
\end{equation}
If the action of \(S_k\) on \(\MULTIPLE{k}{f}\) preserves orientation (i.e.\
\(n-m\) is even or \(k=1\)), then the factor \(\multiple{k}{f}\) is naturally
oriented.

We are going to express the signature and the characteristic numbers of the
manifold of \(k\)-tuple points of \(f\) in terms of cohomological invariants of
\(M\), \(N\) and \(f\) (formulas \eqref{eq:11}, \eqref{eq:13} and
\eqref{eq:17}).

\paragraph{Special cases}
\label{sec:special-cases-1}

Our formulas are particularly simple in many special cases as shown in
Subsection~\ref{sec:special-cases}.  For example, if \(N\) is a sphere, we
recover Sz\H ucs's formulas (here \eqref{eq:71} and \eqref{eq:72}) from
\cite[Theorems~4 and~5]{MR98i:57058}.

%% Hirzebruch's virtual formula
Our formula also generalizes Hirzebruch's virtual signature formula, which we
recall now.

Let \(f\colon M \to N\) be a union of \(k\) embeddings of codimension \(2\) manifolds
\(M_1, \dotsc, M_k\) into \(N\).  Recall that every codimension \(2\) embedded
manifold is the set of zeros of a transversal section of a \(2\) dimensional
vector bundle.  Let \(V_i\) be such a bundle for \(M_i\), and let us denote
the Euler class of this bundle by \(e_i\).  The manifold \(\multiple{k}{f}\)
of \(k\)-tuple points of \(f\) is exactly the intersection of the submanifolds
\(M_1,\dotsc,M_k\).  Hirzebruch's virtual signature formula from
\cite[9.3(4')]{MR34:2573} states that its signature is
\begin{equation}
  \label{eq:3}
  \sigma(\multiple{k}{f}) = \left\langle L(N) \prod_{i=1}^k e_i L(V_i)^{-1}, [N] \right\rangle,
\end{equation}
where \(L(N)\) is the Hirzebruch class of the tangent bundle of \(N\) and,
similarly, \(L(V_i)\) is the Hirzebruch class of \(V_i\).

Hirzebruch uses ``index'' for ``signature'' and also uses a different
notation.  We can get back Hirzebruch's original formula by replacing in the
above formula \(\sigma\) with \(\tau\), the number \(k\) with \(r\), the manifold
\(\multiple{k}{f}\) with \(V^{n-2r}\), the manifold \(N\) with \(M\), the
expression \(\left\langle \mathord{-} , [N] \right\rangle\) with \(\varkappa^n[\mathord{-}]\), the
class \(L(N)\) with \(\sum_{i=0}^\infty L_i(p_1(M^n), \dotsc, p_i(M^n))\), and \(e_i
L(V_i)^{-1}\) with \(\tanh v_i\).  All but the last replacement are just
changes in notation.  To justify the last replacement, note that the dual
class \(v_i\) of \(M_i\) is the Euler class \(e_i\) of \(V_i\) and hence
\begin{equation}
  \label{eq:70}
  e_i L(V_i)^{-1} = e_i \frac{\tanh e_i}{e_i} = \tanh e_i = \tanh v_i.
\end{equation}

Our Corollary~\ref{cor:1} generalizes Hirzebruch's formula to a generic
immersion \(f\colon M \looparrowright N\) with even codimension by introducing the
cohomology class \(B_k(f)\) in~\eqref{eq:19} such that for any immersion
\(f\), we get:
\begin{equation}
  \label{eq:4}
  \sigma(\multiple{k}{f}) = \frac{1}{k!} \langle L(N) B_k(f), [N] \rangle.
\end{equation}
When \(M\) is the disjoint union of manifolds \(M_i\), each of which is
embedded by \(f\), then \(B_k(f)\) factors into a product of some classes of
the \(M_i\) (see Theorem~\ref{th:several-components}).  If all \(M_i\) has
codimension \(2\), then this formula reduces to Hirzebruch's \eqref{eq:3} (see
the discussion after Theorem~\ref{th:several-components}).

Equation~\eqref{eq:4} has a version in \eqref{eq:17} using the cohomology of
\(M\) instead of the cohomology of \(N\).  Sz\H ucs obtained this version in the
special case when \(N\) is a Euclidean space in \cite[Theorem~4]{MR98i:57058},
which we reproduce as \eqref{eq:71}.

\paragraph{The main idea}
\label{sec:main-idea}

Using the definition in~\eqref{eq:1}, let \(i_k\colon \MULTIPLE{k}{f} \to M^{\times k}\)
denote the inclusion and \(j_k\colon \MULTIPLE{k}{f} \to M\) the projection to the
first coordinate.  What we essentially do is deduce an explicit formula for
\(j_{k!}  i_k^*\) (Theorem~\ref{th:technical:jk!ik*}).  The formulas for the
characteristic numbers are applications of this formula.

How do we compute \(j_{k!}  i_k^*\)?  First, in
Section~\ref{sec:from-topol-algebra}, we apply Ronga's clean intersection
theorem \cite[Proposition~2.2]{MR82d:57018} to obtain the
recursion~\eqref{eq:37} for this map.  Then, in
Section~\ref{sec:power-seri-ident}, we interpret this recursion using
generating functions, which produces an easy way to solve it.

\section{Main results}
\label{sec:main-results}

\subsection{Notation}
\label{sec:notation}

We fix our notation for the rest of the paper.
%% Dimensions of manifolds as functions
Let \(f\colon M^{m} \looparrowright N^{n}\) be a generic immersion of compact
oriented smooth manifolds with even codimension.  The dimension of the
components of the manifolds need not be the same.  Let \(m\) be the dimension
function of \(M\) which maps every component of \(M\) to its dimension, and,
similarly, let \(n\) be the dimension function of \(N\).  Now, \(f\) having
even codimension means that for every \(x \in M\) the number \(n(f(x)) -
m(x)\) is even.

\paragraph{Cohomology classes}
\label{sec:cohomology-classes}

We fix notation for some cohomology classes.  Let \(\nu\) be the normal bundle
of \(f\) and let \(e=e(\nu)\) be the Euler class of \(\nu\).  Let \(L(\xi)\) be the
Hirzebruch class of the bundle \(\xi\).  If \(X\) is a manifold, we write \(L(X)\) for
\(L(TX)\), the Hirzebruch class of the tangent bundle of \(X\).  Similarly, we
define \(P(\xi)\) to be the total Pontrjagin class of \(\xi\), and let \(P(X)
\coloneq P(TX)\) be the total Pontrjagin class of the manifold \(X\).  We use
similar notation for Chern classes with \(c\) instead of \(P\).

We need further notation from~\cite{MR98i:57058}.  Let \(J=(j_1, j_2, \dotsc,
j_l)\) be a sequence of non-negative integers.  If \(a\) is a total cohomology
class with \(a_i\) its \(i\) dimensional part, then let \(a_J \coloneq a_{j_1}
a_{j_2} \dotsm a_{j_l}\).  For example, if \(X\) is a compact manifold, then
\(p_J[X] \coloneq \langle P(X)_J, [X] \rangle\) is the Pontrjagin number of \(X\)
corresponding to \(J\) (assuming that \(\sum j_i\) is the dimension of \(X\)).

%% operations among equivalence classes
We write cup products either as ordinary products or with the symbol \(\cup\) for
the operation.  We use \(\times\) for cross products.

\paragraph{Equivalence relations}
\label{sec:equiv-relat}

Equivalence relations naturally arise in our treatment, see
Section~\ref{sec:from-topol-algebra} or \cite{MR98i:57058}, which is the
starting point of our investigation.

Let \(\Eq(k)\) be the set of all equivalence relations on \(\{ 1, \dotsc, k \}\).
We will think of an equivalence relation \(\alpha\) as the set of its equivalence
classes, thus \(\Theta \in \alpha\) will denote that \(\Theta\) is an equivalence class of
\(\alpha\).  Every \(\Theta \in \alpha\) is a subset of \(\{ 1, \dotsc, k\}\) and so \(\Theta\) is
ordered by the usual ordering on integers.  Moreover, there is an ordering
among the equivalence classes themselves: \(\Theta_1 < \Theta_2\) if and only if the
smallest element of \(\Theta_1\) is smaller than the smallest element of \(\Theta_2\).
Hence every equivalence relation \(\alpha\) is the ordered set of its equivalence
classes.

Let \(\alpha[i]\) be the equivalence class of \(\alpha\) containing \(i\).  In
particular, \(\alpha[1]\) is the smallest equivalence class of \(\alpha\).

We will let \(0=0(k)\) denote the trivial equivalence relation, under which
different elements are not equivalent.  Let \(1=1(k)\) denote the universal
equivalence relation, under which all elements are equivalent.

Whenever we write \(\prod_{\Theta \in \alpha}\) we assume that the terms of the product appear
in the order determined by the ordering of \(\alpha\).  Similar remark applies to
\(\prod_{i \in \Theta}\) and other products with ordered index set.

We will denote by \(\size{X}\) the number of elements of the set \(X\).  For
example, \(\size{\alpha}\) is the number of equivalence classes of \(\alpha\).

\paragraph{Maps}
\label{sec:maps}

We list the maps between topological spaces we will use in our formulas.
These are variants of the diagonal map and the graph of \(f^{\times k}\).  Below
\(k\) is a positive integer and \(\alpha\) is an equivalence relation on
\(\{1,\dotsc,k\}\).  Moreover, \(x_{i}\) denotes the \(i\)th coordinate of \(x\).
For an element \(x\) of \(M^{\times \size{\alpha}}\) and \(\Theta \in \alpha\), let \(x_\Theta\) denote
the \(\Theta\)-coordinate of \(x\).  We define \(x_\Theta\) similarly for \(M \times N^{\times
  (\size{\alpha}-1)}\), where \(M\) is the \(\alpha[1]\)-coordinate and the other
coordinates are identified with the other classes of \(\alpha\).

\begin{align}
  \label{eq:5}
  \Delta^{1(k)} = \diagonal{k}\colon M &\to M^{\times k}, &  \diagonal{k}(x)_{i} &\coloneq x\\
  \label{eq:6}
  \Deltaalpha\colon M^{\times \size{\alpha}} &\to M^{\times k}, &  \Deltaalpha(x)_{i} &\coloneq x_{\alpha[i]}\\
  \label{eq:7}
  \Gamma^{1(k)} = \Gamma_k\colon M &\to M \times N^{\times (k-1)}, &  \Gamma_k(x)_{i} &\coloneq 
  \begin{cases}
    x    & \text{if \(i=1\)}\\
    f(x) & \text{if \(i>1\)}
  \end{cases}\\
  \label{eq:8}
  \DMalpha\colon M \times N^{\times (\size{\alpha}-1)} &\to M \times N^{\times (k-1)}, & \DMalpha(x)_{i} &\coloneq
  \begin{cases}
    f(x_{\alpha[i]})  & \text{\pbox[t]{\columnwidth}{if \(\alpha[i]=\alpha[1]\)\\ and \(i>1\)}}\\
    x_{\alpha[i]}     & \text{otherwise}
  \end{cases}
\end{align}
Note that if \(M=N\) and \(f\) is the identity, then \(\DMalpha=\Deltaalpha\).
Occasionally, we will use \(\diagonal{k}\) for manifolds other than \(M\).  The
context will always make this clear.

We need two maps from the multiple point manifold \(\MULTIPLE{k}{f}\) defined
in~\eqref{eq:1}: the canonical inclusion \(i_k\colon \MULTIPLE{k}{f} \to M^{\times k}\)
and the projection \(j_k\colon \MULTIPLE{k}{f} \to M\) to the first coordinate of
\(M^{\times k}\).  When we want to include \(f\) in the notation, we
write \(i_k^{(f)}\) and \(j_k^{(f)}\).

\paragraph{Abbreviation}
\label{sec:abbreviation}

To make formulas more readable, we introduce a shorthand notation for a
frequent constant:
\begin{align}
  \label{eq:69}
  \constant{k} &\coloneq (-1)^{k-1} (k - 1)!
\end{align}

\subsection{The general formula}
\label{sec:general-formula}

Now we state our main results, which will be proved in later sections.

Let \(f\colon M^{m} \looparrowright N^{n}\) be a generic immersion of oriented
compact manifolds with even codimension.
We start with the general formula for signature and characteristic numbers:
\begin{theorem}\label{th:FORMULA}
  The signature and the Pontrjagin numbers of \(\multiple{k}{f}\) are
  \begin{align}
    %% Signature
    \label{eq:9}
      \sigma(\multiple{k}{f}) &= \frac{1}{k!}
      \left\langle
        j_{k!} i_k^* \left(L(M) \times \bigtimes_{i=2}^k L(\nu)^{-1} \right), [M]
      \right\rangle\\
      \label{eq:10}
      &= \frac{1}{k!}
      \left\langle
        L(N) f_!j_{k!} i_k^* \left(\bigtimes_{i=1}^k L(\nu)^{-1} \right), [N]
      \right\rangle,\\
    %% Pontrjagin numbers
    \label{eq:11}
      p_J[\multiple{k}{f}] &= \frac{1}{k!}
      \left\langle
        j_{k!} i_k^* \left(P(M) \times \bigtimes_{i=2}^k P(\nu)^{-1}\right)_J, [M]
      \right\rangle.\\
    %% Chern numbers
\intertext{If \(M\) and \(N\) are almost complex, then we have a similar
  formula for Chern numbers:}
    \label{eq:12}
       c_J[\multiple{k}{f}] &= \frac{1}{k!}
      \left\langle
        j_{k!} i_k^* \left(C(M) \times \bigtimes_{i=2}^k C(\nu)^{-1}\right)_J, [M]
      \right\rangle.
  \end{align}
\end{theorem}
Formulas~\eqref{eq:11} and \eqref{eq:12} can also be written in the form
analogous to~\eqref{eq:10}.

To make these formulas explicit, we have to compute \(j_{k!} i_k^*\) (or \(f_!
j_{k!} i_k^*\)).  Recall that \(e=e(\nu)\) is the Euler class of the normal bundle of \(f\).
\begin{theorem}
  \label{th:technical:jk!ik*}
  For every cohomology class \(x \in H^*(M^{\times k})\)
  \begin{align}
    \label{eq:13}
    j_{k!} (i_k^*(x)) &= \sum_{\alpha \in \Eq(k)} \Gamma_{\size{\alpha}}^* (1 \times f^{\times (\size{\alpha}-1)})_!
    \left(\left(
      \bigtimes_{\Theta \in \alpha} \constant[e]{\size{\Theta}} \right) \cdot \Deltaalpha^* (x)
    \right),\\
    \label{eq:14}
    f_! (j_{k!} (i_k^*(x))) &= \sum_{\alpha \in \Eq(k)} \diagonal{\size{\alpha}}^* f^{\times \size{\alpha}}_!
    \left(\left(
      \bigtimes_{\Theta \in \alpha} \constant[e]{\size{\Theta}}\right) \cdot \Deltaalpha^* (x)
    \right).
  \end{align}
  In particular, if \(x_1, \dotsc, x_k \in H^{2*}(M)\) (i.e.\ the \(x_i\) have
  even dimension), then
  \begin{gather}
    \label{eq:15}
    \begin{split}
      j_{k!} i_k^* (x_1 \times \dotsb \times x_k) = \sum_{\alpha \in \Eq(k)} %\epsilon(\alpha)
      &\left( \constant[e]{\size{\alpha[1]}} \prod_{i \in \alpha[1]} x_i \right)\\
      \cdot \prod_{\substack{\Theta \in \alpha\\ \Theta > \alpha[1]}} &\left( \constant{\size{\Theta}} f^*f_!
        \left( e^{\size{\Theta}-1} \prod_{i \in \Theta}x_i \right)
      \right),
    \end{split}\\
    \label{eq:16}
    f_! j_{k!} i_k^* (x_1 \times \dotsb \times x_k) = \sum_{\alpha \in \Eq(k)} %\epsilon(\alpha)
    \prod_{\Theta \in \alpha} \left( \constant{\size{\Theta}} f_! \left( e^{\size{\Theta}-1} \prod_{i \in \Theta}
        x_i \right) \right).
  \end{gather}
\end{theorem}

The two theorems together provide explicit formulas for the characteristic
numbers and signature.  We state only the signature formula:
\begin{corollary}
  \label{cor:1}
  The signature of the \(k\)-fold intersection manifold of \(f\) is
  \begin{equation}
    \label{eq:17}
    \begin{split}
      \sigma(\multiple{k}{f}) &= \sum_{l + \sum_{i=1}^{k-1} i l_i = k}
      \frac{(-1)^{k-1-\sum_{i=1}^{k-1} l_i}}{k \prod_{i=1}^{k-1} i^{l_i} \cdot l_i!}\\
      &\phantom{=} {} \cdot \left\langle L(M) e^{l -1} L(\nu)^{1-l} \prod_{i=1}^{k-1}
        \left( f^*f_!  \left( e^{i -1} L(\nu)^{-i} \right) \right)^{l_i}, [M]
      \right\rangle\\
      &= \sum_{\sum_{i=1}^k i l_i = k} \frac{(-1)^{k-\sum_{i=1}^k l_i}}{\prod_{i=1}^k i^{l_i} \cdot l_i!}
      \left\langle L(N) \prod_{i=1}^k \left( f_!  \left( e^{i -1} L(\nu)^{-i} \right)\right)^{l_i}
        , [N] \right\rangle
    \end{split}
  \end{equation}
  where the indices \(l_i\) run through the non-negative integers and the
  index \(l\) runs through the positive integers.
\end{corollary}

So the general signature formula is similar to Hirzebruch's
formula~\eqref{eq:3}:
\begin{gather}
  \label{eq:18}
  \sigma(\multiple{k}{f}) = \frac{1}{k!} \langle L(N) B_k(f), [N] \rangle,
  \intertext{where \(B_k(f)\) generalizes the product in~\eqref{eq:3}:}
  \label{eq:19}
  \begin{split}
    B_k(f) &\coloneq f_! j_{k!} i_k^* \left(\bigtimes_{i=1}^k L(\nu)^{-1}\right)\\
    &= \sum_{\sum_{i=1}^k i l_i = k} \frac{k! (-1)^{k-\sum_{i=1}^k l_i}}{\prod_{i=1}^k
      i^{l_i} \cdot l_i!}  \prod_{i=1}^k \left( f_!  \left( e^{i -1} L(\nu)^{-i}
      \right) \right)^{l_i}.
  \end{split}
\end{gather}

%%Several components
We now examine how this formula reduces to Hirzebruch's \eqref{eq:3}, i.e.\ the case
when \(M\) has several components.
\begin{theorem}
  \label{th:several-components}
  Let \(M\) be the disjoint union of manifolds \(M_1, \dotsc, M_l\).  Let
  \(f_i\) denote the restriction of \(f\) to \(M_i\).  Then we can compute
  \(B_k(f)\) as
  \begin{equation}
    \label{eq:20}
    B_k(f) = \sum_{k_1+ \dotsb + k_l = k} \frac{k!}{k_1! \dotsm k_l!} \prod_{i=1}^l B_{k_i}(f_i),
  \end{equation}
  where the \(k_i\) run through the non-noegative integers, and \(B_0(f_i)
  \coloneq 0\) by definition.

  In particular, if the \(M_i\) are embedded manifolds, then
  \begin{equation}
    \label{eq:21}
    B_l(f) = \prod_{i=1}^l B_1(f_i) = \prod_{i=1}^l f_{i!}(L(\nu_i)^{-1})
  \end{equation}
  where \(\nu_i\) is the normal bundle of \(f_i\).
\end{theorem}
Let \(V_i\) be a \(2\) dimensional vector bundle over \(N\).  Let \(M_i\) be
the set of zeros of a transversal section of \(V_i\).  Let \(f_i\) be the
inclusion of \(M_i\) into \(N\).
Then the normal bundle of \(f_i\) is the restriction of \(V_i\) to \(M_i\),
which means \(\nu_i = f_i^*(V_i)\).  Thus we have
\begin{equation}
  \label{eq:68}
  f_{i!}(L(\nu_i)^{-1}) = f_{i!}(f_i^*(L(V_i)^{-1})) = e_i L(V_i)^{-1}
\end{equation}
where \(e_i\) is the Euler class of \(V_i\).  Finally, Hirzebruch's formula
\eqref{eq:3} is obtained by combining \eqref{eq:18}, \eqref{eq:19},
\eqref{eq:21} and \eqref{eq:68}.

%%Proofs
Theorem~\ref{th:technical:jk!ik*} will be proved in
Sections~\ref{sec:from-topol-algebra} and \ref{sec:power-seri-ident}.  The
other results will be proved in Section~\ref{sec:finish-comp}.

\subsection{Special cases}
\label{sec:special-cases}

We present some special cases when the formulas above reduce to a product.
We leave the proofs to Section~\ref{sec:finish-comp}.

\paragraph{\boldmath \protect\(e\protect\), \protect\(L(\nu)\protect\) comes from \protect\(N\protect\)}
\label{sec:e-lnu-comes}

If the cohomology classes \(e\) and \(L(\nu)\) are in the image of \(f^*\), then
the formulas simplify:
\begin{align}
  \label{eq:22}
  j_{k!} i_k^* (f^{\times k})^* (y) &= \diagonal{k}^*(f^{\times k})^* (y) \prod_{i=1}^{k-1} (f^* f_!
  (1) - ie)
  \intertext{leading to}
  \label{eq:23}
  \sigma(\multiple{k}{f}) &= \frac{1}{k!} \left\langle L(M) L(\nu)^{-(k-1)} \prod_{i=1}^{k-1}
    (f^*f_!(1) - ie), [M] \right\rangle\\
  \label{eq:24}
  p_J[\multiple{k}{f}] &= \frac{1}{k!} \left\langle \left( P(M)
      P(\nu)^{-(k-1)} \right)_J \prod_{i=1}^{k-1} (f^*f_!(1) - ie), [M] \right\rangle\\
  \label{eq:25}
  c_J[\multiple{k}{f}] &= \frac{1}{k!} \left\langle \left( C(M)
      C(\nu)^{-(k-1)} \right)_J \prod_{i=1}^{k-1} (f^*f_!(1) - ie), [M] \right\rangle.
\end{align}

\paragraph{\boldmath \protect\(e=0\protect\)}
\label{sec:e=0}

In case \(e=0\), the only non-zero summand in \eqref{eq:13} corresponds to \(\alpha
= 0\).
\begin{align}
  \label{eq:26}
  j_{k!} \circ i_k^* &= \Gamma_k^* \circ (1 \times f^{\times (k-1)})_!
  \intertext{leading to}
  \label{eq:27}
  \sigma(\multiple{k}{f}) &= \frac{1}{k!}
  \left\langle
    L(M) (f^* f_! L(\nu)^{-1})^k, [M]
  \right\rangle\\
  &= \frac{1}{k!} \left\langle L(N) \left( f_!  \left(
        L(\nu)^{-1} \right) \right)^{k}, [N] \right\rangle\\
  \label{eq:28}
  p_J[\multiple{k}{f}] &= \frac{1}{k!} \left\langle \diagonal{k}^* f^{\times k}_! \left(
    P(M) \times \bigtimes_{i=2}^k  P(\nu)^{-1}
    \right)_J , [N] \right\rangle\\
  \label{eq:29}
  c_J[\multiple{k}{f}] &= \frac{1}{k!} \left\langle \diagonal{k}^* f^{\times k}_! \left(
    C(M) \times \bigtimes_{i=2}^k  C(\nu)^{-1}
    \right)_J , [N] \right\rangle.
\end{align}

\paragraph{\boldmath \protect\(f^* f_! = 0\protect\)}
\label{sec:boldm-f-nullh}

If \(f^* f_! = 0\) (this is the case if \(f\) is nullhomotopic and \(n>0\))
then the only non-zero summand in~\eqref{eq:15} corresponds to \(\alpha=1\).  Hence
the formulas reduce to simple products:
\begin{align}
  \label{eq:30}
  j_{k!} i_k^* (x_1 \times \dotsb \times x_k) &= \constant[e]{k} (x_1 \dotsm x_k) =
  (-1)^{k-1} (k-1)! e^{k-1} x_1 \dotsm x_k
  \intertext{leading to}
  \label{eq:31}
  \sigma(\multiple{k}{f}) &= \frac{(-1)^{k-1}}{k} \left\langle e^{k-1} L(M) L(\nu)^{1-k}, [M]
  \right\rangle\\
  \label{eq:32}
  p_J[\multiple{k}{f}] &= \frac{(-1)^{k-1}}{k} \left\langle e^{k-1}
    \left(P(M) P(\nu)^{1-k}\right)_J, [M] \right\rangle\\
  \label{eq:33}
  c_J[\multiple{k}{f}] &= \frac{(-1)^{k-1}}{k} \left\langle e^{k-1}
    \left(C(M) C(\nu)^{1-k}\right)_J, [M] \right\rangle.
\end{align}

In particular, if \(f\) is nullhomotopic, then \(f^*(TN)\) is a trivial bundle
and hence \(L(\nu)^{-1}=L(M)\) and \(P(\nu)^{-1}=P(M)\), so the formulas reduce to
those in \cite[Theorems~4 and~5]{MR98i:57058} (up to minor notational
differences), where these are claimed only when \(N\) is \(\setR^n\) (which we can
replace by the sphere \(S^n\) if we want \(N\) to be compact):
\begin{align}
  \label{eq:71}
  \sigma(\multiple{k}{f}) &= \frac{(-1)^{k-1}}{k} \left\langle e^{k-1} L(M)^{k}, [M]
  \right\rangle\\
  \label{eq:72}
  p_J[\multiple{k}{f}] &= \frac{(-1)^{k-1}}{k} \left\langle e^{k-1}
    \left(P(M)^{k}\right)_J, [M] \right\rangle.
\end{align}
This also corrects a typo (missing sign) in \cite[Theorem~4]{MR98i:57058}.

\section{Sketch of proof}
\label{sec:sketch-proof}

From a technical point of view, the main result is formula~\eqref{eq:13}.  The
other results easily follow from it, as we show in
Section~\ref{sec:finish-comp}.  In this section, we sketch the proof of
formula~\eqref{eq:13}.

Briefly, the proof splits into two parts:
The first part (Section~\ref{sec:from-topol-algebra}) studies the
geometric situation to obtain formula~\eqref{eq:37}.  The second part
is an algebraic rewrite of this formula to achieve our goal: \eqref{eq:13}.

In more detail, as in \cite{MR98i:57058}, the geometric idea is the
description of the preimage of the (narrow) diagonal \(\diagonal{k}(N)\) under \(f^{\times
  k}\).  Its components are parametrized by the equivalence relations on \(k\)
elements.  The component belonging to an equivalence relation with \(l\)
equivalence classes is canonically isomorphic to the manifold of \(l\)-tuple
points of \(f\).  Ronga's \hyperref[th:subcartesian]{Clean Intersection
  Theorem} translates this geometric decomposition into formula~\eqref{eq:37}.

The algebraic manipulation of \eqref{eq:37} is guided by an interpretation of
this formula as a power series equation: \(G = F \circ H\), where \(F\) collects
the unknowns \(j_{k!}  i_k^*\).  The power series \(H\) turns out to be
invertible so we can rewrite the formula as \(F = G \circ H^{-1}\), which is
just~\eqref{eq:13}.

\section{From topology to algebra}% From Sz\H ucs's article
\label{sec:from-topol-algebra}

In this section we derive the recursion~\eqref{eq:37} on \(j_{k!} i_k^*\).

\paragraph{Subcartesian diagram}
\label{sec:subcartesian-diagram}

We will use Ronga's theorem on clean intersections
\cite[Proposition~2.2]{MR82d:57018}.  We recall the notion of clean
intersection:
\begin{definition}\label{def:1}
  Two smooth functions \(f\colon A \to M\) and \(g\colon B \to M\) \emph{intersect cleanly}
  if for every \(a \in A\) and \(b \in B\) such that \(f(a)=g(b)\), there are local
  maps around \(a\) of \(A\), around \(b\) of \(B\) and around \(f(a)=g(b)\)
  of \(M\) such that both \(f\) and \(g\) are linear in these maps.  It
  follows that
  \begin{equation}
    \label{eq:34}
    Z \coloneq \{ (a,b) \in A \times B \mid f(a)=g(b) \}
  \end{equation}
  is a submanifold of \(A \times B\), which we shall call the \emph{clean intersection}
  of \(f\) and \(g\).  The projections of \(Z\) to \(A\) and \(B\)
  form a so called \emph{subcartesian diagram} together with \(f\) and \(g\):
  \begin{diagram}
    \xymatrix{
    %% Spaces
      Z\POS="Z" & A\POS="A"\\    %%% Z -\alpha-> A
      B\POS="B" & M\POS="M"      %%% |      |
    %% Arrows                    %%% \beta      f
      \POS"Z"\ar"A"^{\alpha} \ar"B"^{\beta}   %%% |      |
      \POS"A"\ar"M"^{f}            %%%  V      V
      \POS"B"\ar"M"^{g}            %%% B -g-> M
    }
  \end{diagram}
  The \emph{excess vector bundle} is the bundle \(TM / (TA + TB)\) over \(Z\),
  where we have omitted the obvious pull-back functions in the notation as an
  abuse of language.
\end{definition}

Ronga's theorem states a cohomological identity for subcartesian diagrams:
\begin{theorem}[{Clean Intersection Theorem \cite[Proposition~2.2]{MR82d:57018}}]
  \label{th:subcartesian}
  For every subcartesian diagram we have, using the notation of the above
  definition:
  \begin{equation}
    \label{eq:35}
    g^* (f_!(x)) = \beta_! (e \cdot \alpha^*(x)) \quad (x \in H^*(A)),
  \end{equation}
  where \(e\) is the Euler class of the excess bundle.
\end{theorem}

\paragraph{Main argument}
\label{sec:main-argument}

%% Paragraph summary:
In this paragraph we apply the \hyperref[th:subcartesian]{Clean Intersection
  Theorem} to the maps \(1 \times f^{\times (k-1)}\) and \(\Gamma_k\) to obtain
\eqref{eq:37}.

%% The subcartesian diagram:
%%% The intersection manifold:
First, the clean intersection of the maps is the preimage of the image of
\(\Gamma_k\) under \(1 \times f^{\times (k-1)}\), or, equivalently, the preimage of the
diagonal \(\Delta=\Delta^k(N)\) of \(N^{\times k}\) under \(f^{\times k}\).  As in
\cite{MR98i:57058}, this preimage is the disjoint union of closed submanifolds
\begin{equation}\label{eq:36}
  M_\alpha \coloneq \left\{(x_1, \dotsc, x_k) \in \left( f^{\times k} \right)^{-1} 
  (\Delta) \middle| x_i = x_j \iff i \mathrel{\alpha} j \right\},
\end{equation}
where \(\alpha\) runs over the equivalence relations on \(\{ 1, \dotsc, k \}\).
The manifold \(M_\alpha\) is canonically isomorphic to \(\MULTIPLE{\size{\alpha}}{f}\),
and its inclusion into \(M^{\times k}\) factors as
\begin{equation*}
  \xymatrix@1{M_\alpha \ar[r]^{i_{\size{\alpha}}}& M^{\times \size{\alpha}} \ar[r]^{\Deltaalpha}& M^{\times k}}.
\end{equation*}
Among these, \(M_{0(k)}=\MULTIPLE{k}{f}\) is the manifold whose characteristic
numbers we want to compute.

%%% The maps:
Second, we determine the maps in the subcartesian diagram.  The map from
\(M_\alpha\) into \(M^{\times k}\) is just the canonical embedding.
The map from \(M_\alpha\) to the factor \(M\) of \(\Gamma_k\) is the projection
to  the first coordinate.

%%% The diagram itself:
All in all, the subcartesian diagram of \(1 \times f^{\times (k-1)}\) and \(\Gamma_k\) looks
as below.  The outer square is the subcartesian diagram itself.  The
inner square just explains some maps of the outer square.
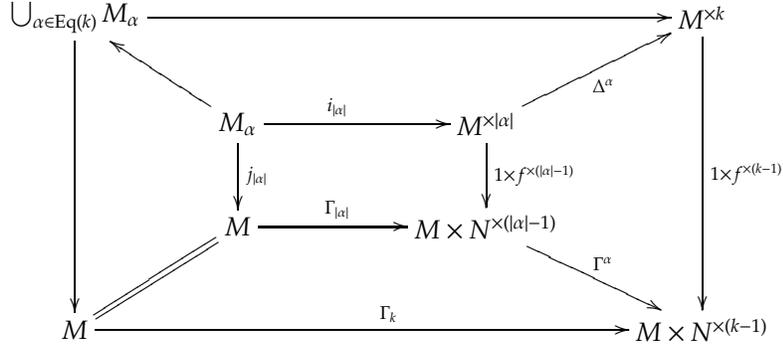
\begin{diagram}
  \xymatrix{
    %% Spaces
    \bigcup_{\alpha \in \Eq(k)} M_\alpha \POS="union" & & & & M^{\times k}        \POS="Mxk"\\
        & M_\alpha   \POS="Ma"  & &  M^{\times \size{\alpha}} \POS="Mxa"\\
        & M     \POS="innerM" & &  M \times N^{\times (\size{\alpha}-1)} \POS="MNxa"\\
    M \POS="outerM" & & & & M \times N^{\times (k-1)}\POS="MNxk"
    %% Arrows
    %%% outer square
    \POS"union"\ar"outerM"
    \POS"outerM"\ar"MNxk"^{\Gamma_k}
    \POS"union"\ar"Mxk"
    \POS"Mxk"\ar"MNxk"^{1 \times f^{\times (k-1)}}
    %%% inner square + arrows to outer
    \POS"Ma"\ar"innerM"^{j_{\size{\alpha}}}        \ar"union"
    \POS"innerM"\ar"MNxa"^(0.4){\Gamma_{\size{\alpha}}}      \ar@2{-}"outerM"
    \POS"Ma"\ar"Mxa"^(0.4){i_{\size{\alpha}}}           %already done: \ar"union"
    \POS"Mxa"\ar"MNxa"^{1 \times f^{\times (\size{\alpha}-1)}} \ar"Mxk"_{\Deltaalpha}
    \POS"MNxa"                       \ar"MNxk"^{\DMalpha}% only to outer
  }
\end{diagram}
%%% The excess vector bundle:
We are going to compute the excess vector bundle.  We will omit the pull-back
maps to simplify our notation since the context will always make it clear
which map is missing.

We fix an equivalence relation \(\alpha\) on \(\{1, \dotsc, k\}\) and determine the
excess vector bundle restricted to \(M_\alpha\).  Therefore we consider all vector
bundles pulled back to \(M_\alpha\).

Recall that the excess vector bundle is the factor of \(T(M \times N^{\times (k-1)})\)
by \(TM^{\times k}\) and \(TM\).  Notice that the inner square of the diagram is a
transverse intersection because \(f\) is generic, so the sum of \(TM^{\times \size{\alpha}}\)
(which is contained in \(TM^{\times k}\)) and \(TM\) is \(T(M \times N^{\times (\size{\alpha}-1)})\).
Thus the excess vector bundle is the factor of \(T(M \times N^{\times (k-1)})\) by
\(TM^{\times k}\) and \(T(M \times N^{\times (\size{\alpha} - 1)})\).

At this point, we notice that the factor makes sense even on \(M^{\times \size{\alpha}}\).
Hence from now on we consider all vector bundles pulled back to \(M^{\times \size{\alpha}}\).

Let \(\Theta_{i}\) denote the \(i\)th equivalence class of \(\alpha\).  We will
write \(l\xi\) for the direct sum of \(l\) copies of a vector bundle \(\xi\).

The embeddings of \(TM^{\times k}\) and \(T(M \times N^{\times (\size{\alpha} - 1)})\) into \(T(M
\times N^{\times (k-1)})\) factor into the components of \(M^{\times \size{\alpha}}\).  For \(i >
1\) on the \(i\)th component we have \(\size{\Theta_{i}}TM\) and \(TN\) embedded
into \(\size{\Theta_{i}}TN\).  The bundle \(TN\) is embedded diagonally, and the
embedding of \(\size{\Theta_{i}}TM\) is induced by \(f\).  Hence the factor is
\((\size{\Theta_{i}}-1)\nu\), where \(\nu\) is the normal bundle of \(f\).  The case of
\(i=1\) is similar.

Thus, the factor on \(M^{\times \size{\alpha}}\) is \((\size{\Theta_{1}}-1)\nu \times\dotsb \times
(\size{\Theta_{l}}-1)\nu\) where \(l \coloneq \size{\alpha}\).  Last, the excess vector
bundle restricted to \(M_\alpha\) is the restriction of this bundle from \(M^{\times
  \size{\alpha}}\) to \(M_\alpha\).

%% The application
Finally, applying the Clean Intersection Theorem
(Theorem~\ref{th:subcartesian}) to the diagram, one obtains:
\begin{equation}
  \label{eq:37}
  \Gamma_k^* (1 \times f^{\times (k-1)})_! (x) = \sum_{\alpha \in \Eq(k)} j_{\size{\alpha}!}
  \left(
    i_{\size{\alpha}}^*
    \left(
      \bigtimes_{\Theta \in \alpha} e^{\size{\Theta}-1}
    \right) \cdot
    i_{\size{\alpha}}^* \Deltaalpha^* x
  \right).
\end{equation}
Recall that \(e=e(\nu)\) is the Euler class of \(\nu\).

\section{The power series identity}
\label{sec:power-seri-ident}

Now we have the recursion formula~\eqref{eq:37} on \(j_{k!} i_k^*\).  To make this
recurrence relation transparent, we interpret it as a power series equality.
Then our theorems will be reduced to routine calculations.

\subsection{Power series}
\label{sec:power-series}

%% The general definition
\paragraph{The general definition}
\label{sec:general-definition}

\emph{Power series} are morphisms of the following category.  Objects
are sequences \(A=(A_k)_{k=1}^\infty\) of modules.  Let \(A=(A_k)_{k=1}^\infty\) and
\(B=(B_k)_{k=1}^\infty\) be two sequences of modules.  A morphism or \emph{power
  series} \(F\) from \(A\) to \(B\) is a collection of homomorphisms \((F_\alpha\colon
A_k \to B_{\size{\alpha}} \mid \alpha \in \Eq(k), k=1,\dotsc,\infty)\).  Given two power series \(F\colon A \to B\)
and \(G\colon B \to C\), we define their \emph{composite} \(G \circ F\) by
\begin{equation}
  \label{eq:38}
  (G \circ F)_\alpha = \sum_{\beta \leq \alpha} G_{\alpha/\beta} \circ F_{\beta}.
\end{equation}
Here \(\beta\) and \(\alpha\) are equivalence relations on the same set.  The notation
\(\beta \leq \alpha\) means that every class of the equivalence relation \(\alpha\) is a union
of some classes of the equivalence relation \(\beta\) (this is the usual ordering
of equivalence relations).  Thus \(\alpha\) induces an equivalence relation \(\alpha/\beta\)
on the classes of \(\beta\): namely, those classes of \(\beta\) are equivalent which
belong to the same class of \(\alpha\).  There is a unique identification between
the ordered set of classes of \(\beta\) and the ordered set \(\{ 1, \dotsc, \size{\beta}
\}\).  Thus we may regard \(\alpha/\beta\) as an equivalence relation on the latter set.
This explains the notation in the above formula.

We leave the easy verifications of the axioms of category to the reader.  The
unit elements are of the form \(E\colon A \to A\) defined as \(E_\alpha = 1\) if \(\alpha = 0 \in
\Eq(k)\) for some \(k\), and \(E_\alpha=0\) for all other \(\alpha\).

%% Example: collection of multilinears from analytics
\paragraph{Classical examples}
\label{sec:classical-examples}

Now we shall see that this definition is an extension of the usual definition
of formal power series.  Classically, given an analytic function \(f\colon U \to V\)
between real vector spaces, its (exponential) power series is the sequence of
its derivatives at a point \(u\) of \(U\).  Let us denote by \((f_k\colon U^{\times k} \to
V : k=1, \dotsc, \infty)\) the \(k\)th derivative of \(f\) at \(u\), it is a
(symmetric) \(k\)-linear map.  In our setting, this corresponds to \(F\colon (U^{\otimes
  k})_{k=1}^\infty \to (V^{\otimes k})_{k=1}^\infty\) defined by
\begin{equation}
  \label{eq:67}
  F_\alpha(u_1 \otimes \dotsb \otimes u_k)
\coloneq \bigotimes_{\Theta \in \alpha} f_{\size{\Theta}}(u_i : i \in \Theta).
\end{equation}
(The arguments of \(f_{\size{\Theta}}\) are the elements \(u_i : i \in \Theta\) in some order.
The order does not matter since the function \(f_{\size{\Theta}}\) is symmetric.)  Our
definition of composition generalizes the composition of usual power series,
since equation~\eqref{eq:38} is the generalization of the formula for the
derivatives of a composite function.

By the above formula~\eqref{eq:67}, we can define for all modules \(U\) and
\(V\) and every sequence \((f_k\colon U^{\times k} \to V : k=1, \dotsc, \infty)\) of multilinear
maps a power series \(F\colon (U^{\otimes k})_{k=1}^\infty \to (V^{\otimes k})_{k=1}^\infty\).  Let us call
power series of this form \emph{classical}.  They are clearly closed under
composition.  We will use classical power series for the \(\setZ[x]\)-modules
\(U=V=\setZ[x]\).  In this special case, every \(k\)-linear function is of the
form \(f_k(x_1, \dotsc, x_k) = a x_1 \dotsm x_k\) for some constant \(a \in
\setZ[x]\).

\subsection{Power series in cohomology}
\label{sec:power-seri-cohom}

%% H^*(M) as a sequence of modules
We are only interested in power series from the sequence of cohomology groups
\((H^*(M^{\times k}))_{k=1}^\infty\) to itself.  We will denote
this sequence by \(\bar{H}^*(M)\).

\paragraph{Special series}
\label{sec:special-series}

We want to map the monoid of classical power series of the \(\setZ[x]\)-module
\(\setZ[x]\) to the monoid of power series of \(\bar{H}^*(M)\).  Substitution of an
element \(e \in H^*(M)\) for \(x\) defines a ring homomorphism \(\setZ[x] \to
H^*(M)\), which maps a polynomial \(a\) into the cohomology class \(a(e)\).
Note that a classical power series of \(\setZ[x]\) is just a sequence
\((a_k)_{k=1}^\infty\) of elements of \(\setZ[x]\).  We map such a sequence to the
power series \(F\) defined by
\begin{equation}
  \label{eq:39}
  F_\alpha(y) \coloneq \left( \bigtimes_{\Theta \in \alpha} a_{\size{\Theta}}(e) \right) \cup \Deltaalpha
  (y)
  \quad (y \in H^*(M^{\times k}), \alpha \in \Eq(k)).
\end{equation}

Note that the exponential power series of many real functions, like the
exponential function \(\exp\) and the natural logarithm \(\ln\), have integer
coefficients and hence are \(\setZ \to \setZ\) and \(\setZ[x] \to \setZ[x]\) series.

%% The power series in the recursion
\paragraph{Solving the recursion}
\label{sec:solving-recursion}

We are now ready to analyze our recursion~\eqref{eq:37}, which we repeat here
in a slightly simpler form:
\begin{equation}
  \label{eq:40}
    \Gamma_k^* (1 \times f^{\times (k-1)})_! (x) = \sum_{\alpha \in \Eq(k)} (j_{\size{\alpha}!} \circ i_{\size{\alpha}}^*)
    \left(
      \bigtimes_{\Theta \in \alpha} e^{\size{\Theta}-1} \cdot \Deltaalpha^* x
    \right).
\end{equation}
Our main observation is that the right-hand side is a special case of the
composition formula~\eqref{eq:38}.  Let us form a power series from the
unknown functions \(j_{k!} \circ i_k^*\):
\begin{equation}
  \label{eq:41}
  F_\alpha \coloneq
  \begin{cases}
    j_{k!} \circ i_k^* & \text{if \(\alpha = 1(k)\) for some \(k\)}\\
    0         & \text{otherwise.}
  \end{cases}
\end{equation}
Note that in the exponential power series expansion of the function
\begin{equation}
  \label{eq:42}
  H(x) \coloneq \frac{\exp(ex)-1}{e} = \sum_{k=1}^\infty \frac{e^{k-1}}{k!} x^k
\end{equation}
the coefficient of the \(k\)th term is \(e^{k-1}\), a polynomial in \(e\).
Therefore we may treat \(H\) as a power series via~\eqref{eq:39} with
\(a_k=e^{k-1}\).  So the right-hand side of~\eqref{eq:40} is just \((F \circ
H)_{1(k)}\).  Clearly, \((F \circ H)_\alpha = 0\) for \(\alpha \neq 1(k)\).  Therefore,
similarly to the definition of \(F\), we can define a power series \(G\) from
the left-hand side of~\eqref{eq:40}:
\begin{equation}
  \label{eq:43}
  G_\alpha \coloneq 
  \begin{cases}
    \Gamma_k^* (1 \times f^{\times (k-1)})_! & \text{if \(\alpha = 1(k)\) for some \(k\)}\\
    0                  & \text{otherwise,}
  \end{cases}
\end{equation}
so that our recursion~\eqref{eq:40} simply means
\begin{equation}
  \label{eq:44}
  G = F \circ H.
\end{equation}

%% Conclusion
The power series \(H\) comes from an invertible function, and hence is
invertible.  The inverse is induced by the inverse function
\begin{equation}
  \label{eq:45}
  H^{-1}(y) = \frac{\ln (1+ey)}{e} = \sum_{k=1}^\infty \frac{(-1)^{k-1}(k-1)!e^{k-1}}{k!}
  y^k = \sum_{k=1}^\infty \frac{\constant[e]{k}}{k!},
\end{equation}
where \(\constant{k} = (-1)^{k-1} (k-1)!\) as declared in \eqref{eq:69}.
We see that the coefficients are polynomials in \(e\) with integer
coefficients and hence we may treat \(H^{-1}\) as a power series from
\(\bar{H}^*(M)\) to itself via~\eqref{eq:39} with \(a_k=\constant[e]{k}\).
Hence the solution of our recursion is \(F = G \circ
H^{-1}\), and this means
\begin{equation}
  \label{eq:46}
  j_{k!} (i_k^*(x)) = \sum_{\alpha \in \Eq(k)} \Gamma_{\size{\alpha}}^* (1 \times f^{\times (\size{\alpha}-1)})_!
    \left(
      \bigtimes_{\Theta \in \alpha} \constant[e]{\size{\Theta}} \cdot \Deltaalpha^* x
    \right).
\end{equation}
This is exactly equation~\eqref{eq:13}.  Applying \(f_!\) to it and using the
identity \(f_! \Gamma_{\size{\alpha}}^* (1 \times f^{\times (\size{\alpha}-1)})_! = \Deltaalpha^*
f_!^{\times \size{\alpha}}\), we obtain \eqref{eq:14}.  Substituting \(x_1 \times \dotsb \times
x_k\) for \(x\) into these two formulas yield \eqref{eq:15} and~\eqref{eq:16}.
We indicate below only how one can deduce \eqref{eq:16} from~\eqref{eq:14}.
Recall that the \(x_i\) are assumed to be even dimensional in
Theorem~\ref{th:technical:jk!ik*}, so no sign appears when we permute them.
\begin{align}
  \label{eq:64}
  \Deltaalpha^*(x_1 \times \dotsb \times x_k) &= %\epsilon(\alpha)
      \bigtimes_{\Theta \in \alpha} \prod_{i \in \Theta} x_i\\
  \label{eq:65}
  \left( \bigtimes_{\Theta \in \alpha} \constant[e]{\size{\Theta}}\right) \left(
    \bigtimes_{\Theta \in \alpha} \prod_{i \in \Theta} x_i \right) &= \bigtimes_{\Theta \in \alpha} \left(
    \constant[e]{\size{\Theta}}\prod_{i \in \Theta} x_i \right)\\
  \label{eq:66}
  f_!^{\times \size{\alpha}} \left( \bigtimes_{\Theta \in \alpha} \constant{\size{\Theta}} e^{\size{\Theta} -
      1}\prod_{i \in \Theta} x_i \right) &= \bigtimes_{\Theta \in \alpha} \constant{\size{\Theta}}
  f_! \left( e^{\size{\Theta} - 1}\prod_{i \in \Theta} x_i \right)
\end{align}

\section{Finishing the computation}
\label{sec:finish-comp}

We have done the hard job in the previous sections.  Now we derive the other
results in Section~\ref{sec:main-results} from
Theorem~\ref{th:technical:jk!ik*} by direct computation.

\subsection{Proof of Theorem~\ref{th:FORMULA}}
\label{sec:proof-theorem-FORMULA}

Recall, e.g.\ from \cite[5.1 and Theorem~8.2.2]{MR34:2573}, that for every
manifold \(X\)
\begin{align}
  \label{eq:47}
  \sigma(X)  &= \langle L(X) , [X] \rangle, \\
  p_J[X] &= \langle P(X)_J, [X] \rangle.
\end{align}

%% The formula with the class [M_{0(k)}] in M^{\times k}

%%% The normal bundle of M_{0(k}} in M^{\times k}
We start by determining the normal bundle of \(\MULTIPLE{k}{f}\) in \(M^{\times k}\)
using the diagram:
%%%% The diagram of transverse preimage
\begin{diagram}
  \xymatrix{
    \MULTIPLE{k}{f} \ar[d]^{p_k} \ar[r]^{i_k} & M^{\times k} \ar[d]^{f^{\times k}}\\
               N    \ar[r]^{\diagonal{k}}            & N^{\times k}
    }
\end{diagram}
where \(p_k \coloneq f \circ j_k\).  Note that \(\MULTIPLE{k}{f}\) is the
transverse preimage of the diagonal of \(N^{\times k}\) under \(f^{\times k}\) at least
in a neighbourhood of \(\MULTIPLE{k}{f}\).  So the normal bundle of
\(\MULTIPLE{k}{f}\) in \(M^{\times k}\) is the pull-back of the normal bundle of
the diagonal in \(N^{\times k}\):
\begin{equation}
  \label{eq:48}
  \nu(i_k) = p_k^* \nu(\Delta) = p_k^* \left( \bigoplus^{k-1} TN \right) = i_k^* ( 1 \times
  \underbrace{f^*TN \times \dotsb  \times f^*TN}_{k-1} ).
\end{equation}

%%% L(\Delta)

Hence one obtains for the Hirzebruch class of \(\MULTIPLE{k}{f}\):
\begin{multline}
  \label{eq:49}
  2 L(\MULTIPLE{k}{f}) = 2 i_k^* (L(M^{\times k})) \cdot L(\nu(i_k))^{-1}\\
 = 2 i_k^* \left( (L(M) \times \dotsb \times L(M)) \cdot (1 \times L(f^*TN)^{-1} \times \dotsb \times L(f^*TN)^{-1})
 \right)\\
 = 2 i_k^* (L(M) \times L(\nu)^{-1} \times \dotsb \times L(\nu)^{-1}).
\end{multline}
We have multiplied everything with \(2\) to get rid of eventual torsion parts.
This has no consequence when computing the signature.
%%% The formula
We get by~\eqref{eq:47}
\begin{equation}
  \label{eq:50}
  \begin{split}
    \sigma(\MULTIPLE{k}{f}) &= \langle i_k^* (L(M) \times L(\nu)^{-1} \times \dotsb \times L(\nu)^{-1}),
    [\MULTIPLE{k}{f}] \rangle\\
    &= \langle j_{k!} i_k^* (L(M) \times L(\nu)^{-1} \times \dotsb \times L(\nu)^{-1}), [M] \rangle.
  \end{split}
\end{equation}
This gives \eqref{eq:9} of Theorem~\ref{th:FORMULA} since \(\sigma(\multiple{k}{f})
= \sigma(\MULTIPLE{k}{f}) / k!\).  Formula \eqref{eq:10} is obtained by using \(2
L(M) = 2 f^*(L(N)) \cdot L(\nu)^{-1}\).  The formulas \eqref{eq:11} and \eqref{eq:12}
are proved similarly.

\subsection{Hirzebruch's virtual signature formula}
\label{sec:hirz-virt-sign}

\paragraph{Proof of Corollary~\ref{cor:1}}
\label{sec:proof-cor:1}

The corollary is obtained by plugging equation~\eqref{eq:15} directly
into~\eqref{eq:9} and plugging \eqref{eq:16} into \eqref{eq:10}.
Substituting \(x_i = L(\nu)^{-1}\) for all \(i\) into~\eqref{eq:16}, the summand
corresponding to \(\alpha\) will depend only on the number of elements of the
classes of \(\alpha\).  Hence we can collect those summands together which are
shown equal by this observation.  Adding the collected terms, we obtain a new
summation whose index will run through all tuples of non-negative integers
\(l_1, \dotsc, l_k\) for which \(\sum_{i=1}^k i l_i =k\), corresponding to the
equivalence relations \(\alpha\) with exactly \(l_i\) pieces of \(i\)-element
classes.  There are exactly \(k!/(\prod_{i=1}^k i!^{l_i} \cdot l_i!)\) many such
equivalence relations.  Hence
\begin{equation}
  \label{eq:51}
  \begin{split}
    f_! j_{k!} i_k^* \left(\bigtimes_{i=1}^k L(\nu)^{-1}\right) &= \sum_{\sum_{i=1}^k i
      l_i = k} \frac{k!}{\prod_{i=1}^k i!^{l_i} \cdot l_i!}  \prod_{i=1}^k \left(
      \constant{i} f_!  \left( e^{i -1} L(\nu)^{-i} \right) \right)^{l_i}\\
    &= \sum_{\sum_{i=1}^k i l_i = k} \frac{k! (-1)^{k-\sum_{i=1}^k l_i}}{\prod_{i=1}^k
      i^{l_i} \cdot l_i!} \prod_{i=1}^k \left( f_! \left( e^{i -1} L(\nu)^{-i} \right)
    \right)^{l_i}
  \end{split}
\end{equation}
Recall from~\eqref{eq:69} that \(\constant{i} = (-1)^{i-1} (i-1)!\), which is
used in the second equation above.

This gives the second formula of Corollary~\ref{cor:1}.  The first formula is
obtained in a similar way using~\eqref{eq:15} but now the equivalence class of
\(1\) is special.  Therefore the summation runs through the tuples \((l, l_1,
\dotsc, l_{k-1})\) corresponding to those equivalence relations for which the
class of \(1\) has \(l\) elements and there are exactly \(l_i\) classes with
\(i\) elements which does not contain \(1\).  The number of such equivalence
relations is \((k-1)!/((l-1)! \prod_{i=1}^k i!^{l_i} \cdot l_i!)\).  Therefore
\begin{equation}
  \label{eq:73}
  \begin{split}
    j_{k!} i_k^*(L(M) \times \bigtimes_{i=2}^k L(\nu)^{-1}) &= \sum_{l + \sum_{i=1}^{k-1} i
      l_i = k} \frac{(k-1)!}{(l-1)! \prod_{i=1}^{k-1} i!^{l_i} \cdot l_i!} \constant[e]{l} L(M)\\
    &\phantom{=}{} \cdot
    L(\nu)^{1-l} \prod_{i=1}^{k-1} \left( \constant{i} f^* f_!
      \left( e^{i-1} L(\nu)^{-i} \right) \right)^{l_i}\\
    &= \sum_{l + \sum_{i=1}^{k-1} i l_i = k}
    \frac{(k-1)! (-1)^{k-1-\sum_{i=1}^k l_i}}{\prod_{i=1}^{k-1} i^{l_i} \cdot l_i!}  L(M)
    e^{l-1}\\
    &\phantom{=}{} \cdot L(\nu)^{1-l} \prod_{i=1}^{k-1} \left( f^* f_!
      \left( e^{i-1} L(\nu)^{-i} \right) \right)^{l_i}.
  \end{split}
\end{equation}

The exponent of \((-1)\) in both formulas of the Corollary is the difference
between \(k\) and the number of equivalence classes.

\paragraph{Proof of Theorem~\ref{th:several-components}}
\label{sec:proof-th:sev-comp}

Now we examine the case when \(M\) is a disjoint union of manifolds
\(M_1,\dotsc,M_l\).  Then the cartesian power \(M^{\times k}\) is the disjoint
union of products \(M_{l_1} \times \dotsb \times M_{l_k}\) for all \(1 \leq l_1, \dotsc,
l_k \leq l\).  This decomposition also decomposes \(\MULTIPLE{k}{f}\), which
leads to a decomposition of \(j_{k!}i_k^*\), and hence \(B_k(f)\), into a sum.
We are gong to determine the summands.

Therefore, let us fix a tuple \((l_1,\dotsc,l_k)\).  Let \(\tilde{j_k}\)
denote the restriction of \(j_k\) to \(\MULTIPLE{k}{f} \cap M_{l_1} \times \dotsb \times
M_{l_k}\), and, similarly, let \(\tilde{i_k}\) denote the restriction of
\(i_k\).

Let \(k_t\) be the multiplicity of \(t\) in the tuple, so that \(M_{l_1} \times
\dotsb \times M_{l_k} \cong \bigtimes_{t=1}^l M_t^{\times k_t}\) by a permutation of
coordinates.  Let \(s\) denote the inclusion of this space into \(M^{\times k}\).

Under this identification, \(\MULTIPLE{k}{f} \cap M_{l_1} \times \dotsb \times M_{l_k}\) is
clearly a subspace of \(\bigtimes_{t=1}^l \MULTIPLE{k_t}{f_t}\).  Actually, it
is the preimage of the diagonal of \(N^{\times l}\) under \(\bigtimes_{t=1}^l (f \circ
j_{k_t}^{(f_t)})\), so the square of the following diagram is a transverse
intersection since \(f\) is generic:
\begin{diagram}
  \xymatrix{
    %% Objects:
    %%% First row
    \MULTIPLE{k}{f} \cap \bigtimes_{t=1}^l M_{l_t}     \POS="cap"&
    \bigtimes_{t=1}^l \MULTIPLE{k_t}{f_t}             \POS="lktft"&&
    \bigtimes_{t=1}^l M_t^{\times k_t}                     \POS="lMtxkt"&
    M^{\times k}                                       \POS="Mxk"\\
    %%% Second row
    N                                          \POS="N"&
    \bigtimes_{t=1}^l N                            \POS="lN"
    %% Arrows:
    %%% From first row:
    \POS"cap"  \ar"lktft"^(0.55){p} \ar"N"^{f \circ \tilde{j_k}}
    \POS"lktft" \ar"lMtxkt"^{\bigtimes_{t=1}^l i^{(f_t)}_{k_t}}
    \ar"lN"^{\bigtimes_{t=1}^l f \circ j^{(f_t)}_{k_t}}
    \POS"lMtxkt" \ar"Mxk"^(0.6){s}
    %%% Inside second row:
    \POS"N" \ar"lN"^{\diagonal{l}}
  }
\end{diagram}
Here \(p\) denotes the obvious inclusion map.

By transversality, we have (e.g. as a special case of
Theorem~\ref{th:subcartesian}):
\begin{equation}
  \label{eq:74}
  f_! \tilde{j_k}_! p^* = \diagonal{l}^* \left( \bigtimes_{t=1}^l f \circ
    j^{(f_t)}_{k_t} \right)_!.
\end{equation}
The composition of the top row is \(\tilde{i_k}\), so we obtain
\begin{equation}
  \label{eq:52}
  f_! \tilde{j_{k}}_{!} \tilde{i_k}^* = f_! \tilde{j_{k}}_{!} p^* \left( \bigtimes_{t=1}^l i^{(f_t)}_{k_t} \right)^*
  s^* = \diagonal{l}^*  \left( \bigtimes_{t=1}^l f \circ j^{(f_t)}_{k_t} \right)_!
  \left( \bigtimes_{t=1}^l i^{(f_t)}_{k_t} \right)^*
  s^*.
\end{equation}
Evaluating the expression at \(\bigtimes_{i=1}^k L(\nu)^{-1}\), we obtain the
summand of \(B_k(f)\) corresponding to \((l_1,\dotsc,l_k)\).  We write \(\nu_t\)
for the normal bundle of \(f_t\), which is the restriction of \(\nu\) to
\(M_t\).
\begin{multline}
  \label{eq:53}
  \diagonal{l}^* \left( \bigtimes_{t=1}^l f \circ j^{(f_t)}_{k_t} \right)_! \left(
    \bigtimes_{t=1}^l i^{(f_t)}_{k_t} \right)^* s^* \left(\bigtimes_{i=1}^k L(\nu)^{-1}\right)\\
  = \diagonal{l}^* \left( \bigtimes_{t=1}^l f \circ j^{(f_t)}_{k_t} \right)_! \left(
    \bigtimes_{t=1}^l i^{(f_t)}_{k_t} \right)^* \left( \bigtimes_{t=1}^l
    (L(\nu_t)^{-1})^{\times k_t} \right)\\ = \prod_{t=1}^l f_{t!} j^{(f_t)}_{k_t!}
  (i^{(f_t)}_{k_t})^* (L(\nu_t)^{-1})^{\times k_t} = \prod_{t=1}^l B_{k_t}(f_t).
\end{multline}
The class \(B_k(f)\) is the sum of the last expression for all tuples
\((l_1,\dotsc,l_k)\), which leads to \eqref{eq:20}.  Note that a tuple
\(k_1,\dotsc,k_l\) appears exactly for \(k!/(k_1! \dotsm k_l!)\) many tuples
\((l_1,\dotsc,l_k)\).

\subsection{Proof of the special cases}
\label{sec:proof-special-cases}

Most of the special cases in Subsection~\ref{sec:special-cases} follow easily
from the general formula, since almost all summands become zero
in~\eqref{eq:13}.  Let us examine equation~\eqref{eq:13} in more detail.  If
\(e=0\) then all summands containing \(e\) become zero and hence only the
summand corresponding to \(\alpha=1\) can be non-zero.  If \(f^* f_! =0\) then we
use~\eqref{eq:15} and see that the only summand which does not contain
\(f^*f_!\) corresponds to \(\alpha=0\), and hence the other summands are zero.

%% Case "Everything from N": not from the general formula
Finally, the case ``\nameref{sec:e-lnu-comes}'' requires more
calculations.  We will use power series again to save some computations.
Let us evaluate \eqref{eq:13} at \(x=1\).
\begin{equation}
  \label{eq:54}
  j_{k!}(1) = \sum_{\alpha \in \Eq(k)} \left(\prod_{\Theta \in \alpha} \constant[e]{\size{\Theta}}
  \right) (f^*f_!(1))^{\size{\alpha}-1} = \partial_3^k q(e, f^*f_!(1), 0).
\end{equation}
where the sum is again a special case of the formula of composition of power
series (or higher order derivatives).  Namely, the two functions we compose
are:
\begin{align}
  \label{eq:55}
  g(y,z) &\coloneq \frac{\exp yz -1}{y}\\
  \label{eq:56}
  h(x,z) &\coloneq \frac{\ln(1+xz)}{x}\\
\intertext{and \(q\) is their composition:}
  \label{eq:57}
  q(x,y,z) &\coloneq  g(y,h(x,z)) =\frac{1}{y}
  \left(
    \exp
    \left(
      y \frac{\ln (1 + xz)}{x}
    \right) -1
  \right).
\end{align}

Now we compute the terms of the power series \(q\) using ordinary power
series.  Recall that
\begin{equation}
  \label{eq:58}
  \exp(t \ln(1+x)) = (1+x)^t = \sum_{n=0}^\infty \binom{t}{n} x^n.
\end{equation}
Substituting \(y/x\) for \(t\) and \(xz\) for \(x\) this becomes
\begin{equation}
  \label{eq:59}
  1 + y q(x,y,z) = \exp \left( y \frac{\ln (1+xz)}{x} \right) = \sum_{n=0}^\infty \binom{y/x}{n} x^nz^n.
\end{equation}
Thus the \(n\)th partial derivative of \(q\) in its third variable \(z\) is
\begin{equation}
  \label{eq:60}
  \partial_3^n q(x,y,0) = \frac{n!}{y} \binom{y/x}{n} x^n = \prod_{i=1}^{n-1} (y - ix).
\end{equation}

Plugging this into~\eqref{eq:54}, we obtain
\begin{gather}
  \label{eq:61}
  j_{k!}(1) = \prod_{i=1}^{k-1} (f^*f_!(1)-ie),\\
  \label{eq:62}
  j_{k!} i_k^* (f^{\times k})^* (y) =  \diagonal{k}^*(f^{\times k})^* (y) j_{k!}(1) = \diagonal{k}^*(f^{\times k})^* (y) \prod_{i=1}^{k-1} (f^* f_! (1) -ie).
\end{gather}

The last formula is exactly \eqref{eq:22}, from which \eqref{eq:23},
\eqref{eq:24} and~\eqref{eq:25} are straightforward.

\bibliography{multiple}\bibliographystyle{amsplain}
\end{document}